\documentclass[11pt]{article}    
\usepackage{cite}
\usepackage[pdftex]{graphicx}
\usepackage[cmex10]{amsmath}
\usepackage[tight,footnotesize]{subfigure}
\usepackage{amsfonts,latexsym,amssymb,amsmath}
\usepackage{url,epstopdf}
\newtheorem{algo}{Algorithm}
\newtheorem{proposition}{Proposition}
\newtheorem{corol}{Corollary}
\selectfont
\begin{document}
\title{Some Special Cases in the Stability Analysis of Multi-Dimensional Time-Delay Systems Using The Matrix Lambert W function}
\author{Rudy Cepeda-Gomez\thanks{Institute of Automation, University of Rostock, Rostock, Germany. \texttt{rudy.cepeda-gomez@uni-rostock.de}} and Wim Michiels\thanks{Department of Computer Science, KU Leuven, Leuven, Belgium. \texttt{wim.michiels@cs.kuleuven.be}}}
\maketitle                           

\begin{abstract}
This paper revisits a recently developed methodology based on the matrix Lambert W function for the stability analysis of linear time invariant, time delay systems. By studying a particular, yet  common, second order system, we show that in general there is no one to one correspondence between the branches of the matrix Lambert W function and the characteristic roots of the system. Furthermore, it is shown that under mild conditions only two branches suffice to find the complete spectrum of the system, and that the principal branch can be used to find several roots, and not the dominant root only, as stated in previous works. The results are first presented analytically, and then verified by numerical experiments.
\end{abstract}
\section{Introduction}\label{intro}
In this paper we consider Linear Time Invariant-Time Delay Systems (\textsc{lti-tds}), represented by Delay-Differential Equations (\textsc{ddes}) of the form:
\begin{equation}
\dot{\mathbf{x}}\left(t\right)=\mathbf{Ax}\left(t\right)+\mathbf{Bx}\left(t-\tau\right)
\label{eq:lti-tds}
\end{equation}
The stability analysis and control synthesis of this class of systems is a wide open area of research. The difficulty of this problem arises from the the fact that the delay makes this class of systems infinite dimensional. A nice review of the recent results and challenges int his area can be found in \cite{Sipahi2011}. Several avenues have been followed by different researchers to find solutions to the stability and control questions. Some works study the absolute stability regions with respect to the time delay \cite{Olgac2002,Sipahi2006}, and lead to control strategies that use the time delay as a stabilizing tool \cite{Olgac2005}. Other researchers focus on the numerical computation of the characteristic roots of the system. These works include approaches  based on the discretization of the solution operator\cite{Engelborghs2002,Breda2006,Butcher2011} or its infinitesimal generator \cite{Wu2012,Breda2009}, and methods based on root finding of the characteristic equation \cite{Vhylidal2009}. From here methods have been proposed  to optimize the location of the dominant roots to guarantee a certain performance \cite{Michiels2002,Mondie2003}. Finally, a Krylov method for computing characteristic roots of large-scale problems has been proposed in \cite{Jarlebring2010}.

In the past decade a framework for analyzing \textsc{ddes} based on the Lambert W function has been developed \cite{Asl2003,Yi2009,Yi2010d}. It expands the earlier work \cite{Wright1959}. The main idea of the methodology is to express the solution of a \textsc{dde} as the sum of a series of infinitely many exponential functions. The characteristic roots of the system are found analytically in terms of the Lambert W function. While the problem remains infinite dimensional, a one-to-one correspondence between the characteristic roots of the system and the branches of this multi-valued function is assumed. The stability question is then solved by earmarking the dominant characteristic roots of the system with the branches of the Lambert W function corresponding to $k=0,\,\pm1,\,\ldots,\pm m$, where $m$ is the nullity of the matrix $\mathbf{B}$ in \eqref{eq:lti-tds}. Therefore, only a few branches have to be considered to determine whether a solution is stable or not. Furthermore, the existence of an explicit solution expressed in terms of a power series, allows analyzing structural properties like observability and controllability for systems of \textsc{ddes} \cite{Yi2008}, the development of pole placement techniques for control synthesis \cite{Yi2010a,Yi2010b,Yi2010c,Yi2012b} and other applications like the estimation of decay rates \cite{Duan2012} and spectrum design \cite{Wei2014}.

The basic foundation of his methodology, i.e., the assumption that the principal branch of the Lambert W function defines the stability of the system, is well established for first order systems \cite{Shinozaki2006,Asl2003}. For higher order systems, however, this result has not been extended with the same rigor, and is rather based upon observations \cite{Yi2006a,Yi2006b}.  

In this paper, we show that the main assumption does not hold in general. By studying  particular second order delay systems, with a structure that is very common in applications, we show that the full spectrum of a \textsc{lti-tds} of this class can be found using only two branches of the Lambert W function, i.e., there is no one-to-one correspondence between the characteristic roots and the branches of the Lambert W function. This is due to the fact that an important nonlinear equation in the approach does not have a unique solution. Furthermore, we show that the principal branch can be used to find not only the dominant root of the system, but some other roots too.

The organization of the paper is as follows. In section \ref{sec:LambW} we present the Lambert W function and its matrix version, required to work with higher order systems. Section \ref{sec:solLamb} reviews the methodology to solve \textsc{ddes} using the Lambert W function, both in the scalar and vector cases. Section \ref{sec:args} presents the analysis of a second order system and contains the main results. This analysis is illustrated by numerical examples in Section \ref{sec:numer}. Section \ref{sec:odd} presents discussions on another special case for which the results of section \ref{sec:args} cannot be directly applied. Finally, some conclusions of the study are given in section \ref{sec:conclusions}. 

In the remainder of the paper, scalar quantities are denoted by italic symbols ($a,\,b,\,\lambda$) whereas vectors and matrices are represented by bold face lowercase ($\mathbf{x}$) and uppercase ($\mathbf{A},\,\mathbf{B}$) letters, respectively. The notation $e^{b}$ is used to represent the exponential function of a scalar and $\exp\left(\mathbf{A}\right)$ represents a matrix exponential function. In a similar way, $W_k\left(z\right)$ represents the $k$-th branch of the Lambert W function of a scalar number, and $\mathbf{W}_k\left(\mathbf{A}\right)$ is the matrix Lambert W function.

\section{The Lambert W Function}
\label{sec:LambW}
\begin{figure}
\centering
\includegraphics[scale=0.55]{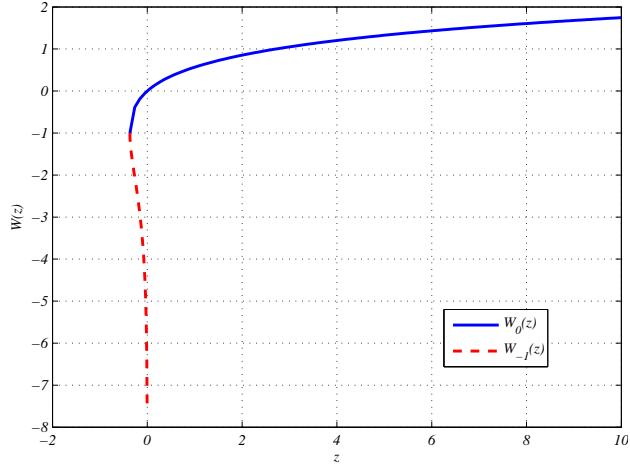}
\caption{The two real branches of the Lambert W function.}
\label{fig:realbranches}
\end{figure}

The Lambert W function is a function $W(z)$, $\mathbb{C}\mapsto\mathbb{C}$, defined as the solution  to the equation
\begin{equation}
W\left(z\right)e^{W\left(z\right)}=z
\label{eq:lambert_scalar}
\end{equation}
This is a multi-valued function, that is, for a $z\in\mathbb{C}$ there are infinitely many solutions to \eqref{eq:lambert_scalar}. To identify these values a branch number is assigned, and we refer to $W_k\left(z\right)$ as the $k$-th branch of the Lambert W function of $z$. The branch cuts are defined in such a way that each branch has a precisely defined range \cite{Corless1996}. For $z\in\mathbb{R}$, only two of the branches are real valued. The principal branch, $W_0\left(z\right)$ is real for $z\geq -1/e$ and its range is the interval $\left[-1,\,\infty\right)$. The branch $W_{-1}\left(z\right)$ is real for $-1/e\leq z<0$, and its range is $\left(-\infty,\,-1\right]$. These branches are shown in Figure \ref{fig:realbranches}.

A comprehensive study of the history, definition and properties of the Lambert W function is found in \cite{Corless1996}. 

The matrix Lambert W function is defined now. Consider a matrix $\mathbf{H}\in\mathbb{C}^{n\times n}$, which has the Jordan canonical decomposition $\mathbf{H}=\mathbf{ZJZ}^{-1}$, with $\mathbf{J}=\mathrm{diag}\left(\mathbf{J}_1\left(\lambda_1\right),\,\mathbf{J}_2\left(\lambda_2\right),\,\ldots,\,\mathbf{J}_p\left(\lambda_p\right)\right)$. Following one of the standard definitions for a function of a matrix \cite{Higham2008}, the Matrix Lambert W function for a Jordan block of size $m$ is defined as:
\begin{equation}
\begin{split}
&\mathbf{W}_k\left(\mathbf{J}_i\right)=\\
&\left[\begin{array}{cccc}
W_k\left(\lambda_i\right)&W_k'\left(\lambda_i\right)&\cdots&\frac{1}{\left(m-1\right)!}W_{k}^{m-1}\left(\lambda_i\right)\\
0&W_k\left(\lambda_i\right)&\cdots&\frac{1}{\left(m-2\right)!}W_{k}^{m-2}\left(\lambda_i\right)\\
\vdots&\vdots&\ddots&\vdots\\
0&0&\cdots&W_k\left(\lambda_i\right)\end{array}\right]
\end{split}
\label{eq:wjordan}
\end{equation}
and the matrix Lambert W function of $\mathbf{H}$ is defined as:
\begin{equation}
\begin{split}
\mathbf{W}_k\left(\mathbf{H}\right)=&\mathbf{Z}\,\mathrm{diag}\left(\mathbf{W}_k\left(\mathbf{J}_1\left(\lambda_1\right)\right),\mathbf{W}_k\left(\mathbf{J}_2\left(\lambda_2\right)\right),\ldots,\right.\\
&\left.\mathbf{W}_k\left(\mathbf{J}_p\left(\lambda_p\right)\right)\right)\mathbf{Z}^{-1}.
\end{split}
\label{eq:wmatrix}
\end{equation}

Every matrix defined by \eqref{eq:wmatrix}, for $k=0,\,\pm 1,\,\pm 2\,\ldots$, is one particular solution to the matrix equation
\begin{equation}
\mathbf{W}_k\left(\mathbf{H}\right)\exp\left(\mathbf{W}_k\left(\mathbf{H}\right)\right)=\mathbf{H}.
\label{eq:wlmatrix}
\end{equation}
The above, standard definition implies that the same branch of the Lambert W function is used in each Jordan block. This is not necessary to have a solution of \eqref{eq:wlmatrix}. Since $W_0\left(0\right)=0$ and $W_k\left(0\right)=\infty$ for $k\neq 0$, we deviate from the standard definition to avoid the infinite value. This special case, called the \emph{hybrid branch} case, is defined in \cite{Yi2009,Yi2010d}. More precisely, there $\mathbf{W}_k\left(\mathbf{H}\right)$ uses the actual value of $k$ for those Jordan blocks with $\lambda\neq 0$ and $k=0$ for those blocks in which $\lambda=0$. This definition is used in the remainder of the paper. 

In a similar way, it is assumed that $e^{-1}$ is not eigenvalue of $\mathbf{H}$ corresponding to a Jordan block of dimension higher than 1 when the principal branch is being computed. This is required to overcome the difficulty represented by the fact that $W_0'\left(e^{-1}\right)$ is not defined. This limitation reduces the elegance of the definition of the matrix Lambert W function, but does not affect its usage \cite{Jarlebring2007}.
\section{Solution of Delay-Differential Equations using the Lambert W Function}
\label{sec:solLamb}
\subsection{Scalar Case}
For a scalar, homogeneous \textsc{dde},
\begin{equation}
\dot{x}\left(t\right)=ax\left(t\right)+b\left(t-\tau\right),
\label{eq:scalardde}
\end{equation}
the characteristic equation is
\begin{equation}
s-a-be^{-s\tau}=0.
\label{eq:cescalar}
\end{equation}
The solution to \eqref{eq:cescalar} can be expressed in terms of the Lambert W function following simple steps \cite{Corless1996,Asl2003,Yi2009,Yi2010d}. This solution has the form
\begin{equation}
s_k=\frac{1}{\tau}W_k\left(\tau be^{-a\tau}\right)+a,
\label{eq:solscalar}
\end{equation} 
where $k=0,\,\pm 1,\,\pm 2,\,\ldots$ indicates the branch of the Lambert W function to be used. Each one of the infinitely many roots of \eqref{eq:cescalar} corresponds to one of the branches of this function. 

It has been proven \cite{Shinozaki2006} that among all  solutions in \eqref{eq:solscalar}, the one that corresponds to the principal branch, $k=0$, always has the largest real part and, therefore, introduces the dominant mode to the solution of the equation. To study the stability of the solution to a one dimensional \textsc{dde} as \eqref{eq:scalardde}, it is necessary and sufficient to find only the solution of \eqref{eq:cescalar} in \eqref{eq:solscalar} corresponding to $k=0$.   
\subsection{Higher Order case}
Consider now a higher order \textsc{DDE} described by:
\begin{equation}
\dot{\mathbf{x}}\left(t\right)=\mathbf{Ax}\left(t\right)+\mathbf{Bx}\left(t-\tau\right)
\label{eq:ltitds}
\end{equation}
with $\mathbf{x}\left(t\right)\in\mathbb{R}^{n}$, $\mathbf{A},\,\mathbf{B}\in\mathbb{R}^{n\times n}$ and $\tau>0$.

The following steps, introduced in \cite{Asl2003} and extended in \cite{Yi2009,Yi2010d,Yi2006a,Yi2006b}, aim at computing characteristic roots using the matrix Lambert W function. The proposed method is based on finding a solution of the equation
\begin{equation}
\mathbf{S}-\mathbf{A}-\mathbf{B}\exp\left(-\mathbf{S}\tau\right)=\mathbf{0},
\label{eq:char}
\end{equation}
where $\mathbf{S}\in\mathbb{C}^{n\times n}$. A matrix $\mathbf{Q}\in\mathbb{C}^{n\times n}$ is introduced, such that
\begin{equation}
\tau\left(\mathbf{S}-\mathbf{A}\right)\exp\left(\left(\mathbf{S}-\mathbf{A}\right)\tau\right)=\tau\mathbf{B}\mathbf{Q}
\label{eq:4}
\end{equation}
is satisfied. Let us define $\mathbf{M}:=\tau\mathbf{B}\mathbf{Q}$. Then, from \eqref{eq:4} and \eqref{eq:wlmatrix}, 
\begin{equation}
\mathbf{S}_k=\frac{1}{\tau}\mathbf{W}_k\left(\mathbf{M}\right)+\mathbf{A},
\label{eq:5}
\end{equation}
with $k\in\mathbb{Z}$, is a solution of \eqref{eq:4}. By substituting \eqref{eq:5} into \eqref{eq:char} the following expression is obtained,
\begin{equation}
\mathbf{W}_k\left(\mathbf{M}\right)\exp\left(\mathbf{W}_k\left(\mathbf{M}\right)+\mathbf{A}\tau\right)-\tau\mathbf{B}=\mathbf{0}.
\end{equation}
Therefore, the steps to compute characteristic roots are given by the following algorithm \cite{Yi2009,Yi2010d}.
\begin{algo}\label{algo}
Repeat for $k=0,\,\pm 1,\,\pm 2,\,\ldots$:
\begin{enumerate}
\item Solve the nonlinear equation
\begin{equation}
\mathbf{W}_k\left(\mathbf{M}_k\right)\exp\left(\mathbf{W}_k\left(\mathbf{M}_k\right)+\mathbf{A}\tau\right)-\tau\mathbf{B}=\mathbf{0},
\label{eq:findQ}
\end{equation}
for $\mathbf{M}_k$\ ( $=\tau \mathbf{B}\mathbf{Q}_k$).
\item Compute $\mathbf{S}_k$ corresponding to $\mathbf{M}_k$ as 
\begin{equation}
\mathbf{S}_k=\frac{1}{\tau} \mathbf{W}_k(\mathbf{M}_k)+\mathbf{A}.
\end{equation}
\item Compute the eigenvalues of $\mathbf{S}_k$.
\end{enumerate} 
\end{algo}

For system \eqref{eq:ltitds} to be stable, all characteristic roots  must have negative real parts. Calculating the solution for all the branches is not possible. To work around this difficulty, the proposers of this methodology assume that for any branch $k$ of the matrix Lambert W function, there is a unique solution $\mathbf{M}_k$ to \eqref{eq:findQ} and a corresponding $\mathbf{S}_k$ matrix. This assumption, based on observations from many examples, leads to a stronger conjecture: when the rank of $\mathbf{B}$ in \eqref{eq:lti-tds} is at least $n-1$, i.e., $\mathbf{B}$ does not have a repeated zero eigenvalue, the characteristic roots with largest real part correspond to the $\mathbf{S}_0$ matrix, found using the principal branch of the matrix Lambert W function in Algorithm~\ref{algo}. This conjecture is formally stated in \cite{Yi2009} and it is the basis for several derivative works \cite{Yi2010a,Yi2010b,Yi2010c,Yi2012a,Yi2012b,Duan2012,Wei2014}.

In the following section, we show that this conjecture does not hold in general. By studying a particular, yet very common, control system, we can show that (\ref{eq:findQ}) may have multiple solutions and there is no one to one correspondence between the branches of the matrix Lambert W function and the characteristic roots of the system. Furthermore, we show that in particular cases it is possible to find \emph{all} characteristic roots using only \emph{two} branches, those corresponding to $k=0$ and $k=-1$. This contradicts the conjecture regarding stability.  

\section{A Common Special Case in Second Order Systems}
\label{sec:args}
We consider a single input, second order system under time delayed state feedback. This type of control system is ubiquitous in applications. Without loss of generality, we assume that the system matrix is in companion form. This system is represented by \eqref{eq:ltitds} with the following $\mathbf{A}$ and $\mathbf{B}$ matrices: 
\begin{equation}
\mathbf{A}=\left[\begin{array}{cc}0&1\\a_{21}&a_{22}\end{array}\right]\quad\mathbf{B}=\left[\begin{array}{cc}0&0\\b_{21}&b_{22}\end{array}\right].
\label{eq:sys2nd}
\end{equation}
Since $\mathbf{B}$ in \eqref{eq:sys2nd} has nullity 1, the current theory predicts that the solutions to \eqref{eq:findQ} corresponding to $k=0$ and $k=\pm1$ generate only its dominant roots.

From the structure of $\mathbf{B}$, we can see that $\mathbf{M}_k=\tau\mathbf{B}\mathbf{Q}_k$, for any given $\mathbf{Q}_k$, has the form:
\begin{equation}
\mathbf{M}_k=\left[\begin{array}{cc}0&0\\m_{21}&m_{22}\end{array}\right].
\label{eq:m}
\end{equation}
Applying the definition of the matrix Lambert W function, using the hybrid branch case because $\mathbf{M}_k$  has one eigenvalue equal to zero, we get for $m_{22}\neq 0$:
\begin{equation}
\mathbf{W}_k\left(\mathbf{M}_k\right)=\left[\begin{array}{cc}0&0\\\frac{m_{21}}{m_{22}}W_k\left(m_{22}\right)&W_k\left(m_{22}\right)\end{array}\right],
\label{eq:wm}
\end{equation}
and from here we obtain
\begin{equation}
\begin{split}
\mathbf{S}_k&=\frac{1}{\tau}\mathbf{W}_k\left(\mathbf{M}_k\right)+\mathbf{A}\\
&=\left[\begin{array}{cc}0&1\\\frac{m_{21}}{\tau m_{22}}W_k\left(m_{22}\right)+a_{21}&\frac{1}{\tau}W_k\left(m_{22}\right)+a_{22}\end{array}\right].
\end{split}
\label{eq:s}
\end{equation}
In case $m_{22}=0,\ m_{21}\neq 0$, a simple computation yields
\begin{equation}
\mathbf{S}_k=\left[\begin{array}{cc} 0 &1\\ \frac{m_{21}}{\tau}+a_{21} & a_{22}\end{array}\right],
\end{equation}
where we used $W_0^{\prime}(0)=1$. We are now ready to state the main results of the paper.

\begin{proposition} \label{prop1}
Let $A$ and $B$ be given by (\ref{eq:sys2nd}). Let $\{\lambda,\,\bar{\lambda}\}$ be any pair of complex conjugate characteristic roots of the system defined by \eqref{eq:sys2nd}. Assume their multiplicity is one. Then for either $k=0$ or $k=-1$ there exists a real solution of (\ref{eq:findQ}), such that, if this solution and corresponding value of $k$ are selected in the first step of Algorithm~\ref{algo}, the characteristic roots $\lambda$ and $\bar{\lambda}$ are found in the last step of the algorithm.
\end{proposition}

\noindent \textbf{Proof.} The key idea is to perform the steps of Algorithm~\ref{algo} in \emph{reverse} order, in the course of which $k$ is selected.

From the pair $(\lambda,\bar{\lambda})$ we first construct a real matrix $\mathbf{S}_k$, of which they are the eigenvalues, namely
\begin{equation}
\mathbf{S}_k=\left[\begin{array}{cc}0&1\\-\left|\lambda\right|^2&2\Re\left(\lambda\right)\end{array}\right].
\label{eq:skproof} 
\end{equation}

Subsequently, we construct $\mathbf{M}_k$ from (\ref{eq:skproof}), where we make distinction between two cases.

\noindent \emph{Case 1: $2\Re(\lambda)\neq a_{22}$.}
Comparing \eqref{eq:s} and \eqref{eq:skproof}, we can take $\mathbf{M}_k$ of the form (\ref{eq:m}), where $m_{21}\in\mathbb{R}$ and $m_{22}\in\mathbb{R}$ are chosen such that the following equations are satisfied:
\begin{align}
W_k\left(m_{22}\right)&=\tau\left(2\Re\left(\lambda\right)-a_{22}\right),\label{eq:wk}\\
m_{21}&=-\frac{m_{22}\left(\left|\lambda\right|^2+a_{21}\right)}{2\Re\left(\lambda\right)-a_{22}}.
\end{align}
Such a choice  is always possible for either $k=0$ or $k=-1$, which then fix $k$.  Equation \eqref{eq:wk} namely implies that $W_k\left(m_{22}\right)$ must be a real number. As mentioned earlier, the definition of the branch cuts of the Lambert W function makes this function to have two real branches: $k=0$, the principal branch, and $k=-1$, and the union of the ranges of $W_0$ and $W_{-1}$ for real valued arguments includes $\mathbb{R}$, see Figure~\ref{fig:realbranches}. Furthermore, in the case considered, $m_{22}$ computed from (\ref{eq:wk}) is different from zero, justifying the use of (\ref{eq:s}).

\noindent \emph{Case 2: $2\Re(\lambda)=a_{22}$.} We can freely choose $k\in\{0,-1\}$ and take 
\begin{equation}
\mathbf{M}_k=\left[\begin{array}{cc}0 & 0 \\ -\tau(a_{21}+|\lambda|^2) & 0\end{array}\right].
\end{equation}

Finally, it remains to show that the constructed pair $(k,\mathbf{M}_k)$ is a solution of (\ref{eq:findQ}). Let $\mathbf{v}$, respectively $\bar{\mathbf{v}}$, be the eigenvector corresponding to $\lambda$, respectively $\bar\lambda$. If the characteristic root $\lambda$ is simple, then the pair $(\mathbf{V},\boldsymbol{\Lambda})$ is an invariant pair of (\ref{eq:ltitds}), where
\begin{equation}
\mathbf{V}=[\mathbf{v}\ \bar{\mathbf{v}}],\ \ \boldsymbol{\Lambda}=\mathrm{diag}(\lambda,\bar\lambda).
\end{equation}
As a consequence it satisfies
\begin{equation}\label{invariant1}
\mathbf{V} \boldsymbol{\Lambda}-\mathbf{A}-\mathbf{B} \mathbf{V} \exp(-\boldsymbol{\Lambda})=0,
\end{equation}
see \cite{Beyn2011}.
The eigenvalue decomposition of $\mathbf{S}_k$ in (\ref{eq:skproof}) takes the form
\begin{equation}
\mathbf{S}_k=\mathbf{V} \boldsymbol{\Lambda} \mathbf{V}^{-1},
\end{equation}
from which we have
\begin{equation}
\exp(-\mathbf{S}_k)=\mathbf{V} \exp(-\boldsymbol{\Lambda}) \mathbf{V}^{-1}.
\end{equation}
It follows that
\begin{equation}
\mathbf{V}\boldsymbol{\Lambda}= \mathbf{S}_k \mathbf{V},\ \ \mathbf{V} \exp(-\boldsymbol{\Lambda})=\exp(-\mathbf{S}_k) \mathbf{V}.
\end{equation} 
Substituting the latter in (\ref{invariant1}) yields, as $V$ is invertible,
\begin{equation}
\mathbf{S}_k-\mathbf{A}-\mathbf{B}\exp(-\mathbf{S}_k)=0.
\end{equation}
Finally, replacing $\mathbf{S}_k$ by
\begin{equation}
\frac{1}{\tau}\mathbf{W}_k\left(\mathbf{M}_k\right)+\mathbf{A}
\end{equation}
results in (\ref{eq:findQ}).
\hfill $\Box$

\begin{proposition} \label{prop2}
Let $A$ and $B$ be given by (\ref{eq:sys2nd}). Let $\lambda_1\neq \lambda_2$ be two real, simple characteristic roots of the system defined by \eqref{eq:sys2nd}. Then for either $k=0$ or $k=-1$ there exists a real solution of (\ref{eq:findQ}), such that, if this solution and corresponding value of $k$ are selected in the first step of Algorithm~\ref{algo}, the characteristic roots $\lambda_1$ and $\lambda_2$ are found in the last step of the algorithm.
\end{proposition}

\noindent\textbf{Proof.}\  The proof is completely analogous to the proof of Proposition~\ref{prop1}. The differences are that we start by defining
\begin{equation}
\mathbf{S}_k=\left[\begin{array}{cc} 0 &1 \\-\lambda_1\lambda_2 & \lambda_1+\lambda_2 \end{array}\right],
\end{equation}
and that the two cases to be considered are $\lambda_1+\lambda_2=a_{22}$ and $\lambda_1+\lambda_2\neq a_{22}$. \hfill $\Box$

From Propositions \ref{prop1}-\ref{prop2} the following Corollary can be derived.
\begin{corol}\label{cor1}
Let $A$ and $B$ be given by (\ref{eq:sys2nd}).  If all characteristic roots of \eqref{eq:sys2nd} are simple and if the number of real characteristic roots  is different from one, then \emph{all} characteristic roots can be found using only two branches of the matrix Lambert W function in Algorithm~\ref{algo}, namely $k=0$ and $k=-1$. Moreover, one can restrict to the real solutions of (\ref{eq:findQ}).
\end{corol}

\textbf{Remark} The previous analysis shows that, with proper initial conditions, all the characteristic roots of system l\eqref{eq:sys2nd} can be found using the branches corresponding to $k=0$ and $k=-1$. However, higher branch numbers can be used to find different pairs of roots following the same reverse engineering approach, without any particular structure. This is demonstrated in the following section. 

\section{Numerical Examples}
\label{sec:numer}
We consider a system defined by matrices:
\begin{equation}
\mathbf{A}=\left[\begin{array}{rr}0&1\\-5&-1\end{array}\right]\quad\mathbf{B}=\left[\begin{array}{rr}0&0\\-3&-0.6\end{array}\right],
\label{eq:examp}
\end{equation}
with a value $\tau=5$ for the time delay.

When the \emph{LambertDDE} toolbox \cite{Yi2012a} was used to calculate the characteristic roots of the system, the software was not able to find a solution of (\ref{eq:findQ}) for any value of $k$ with the default settings. For the numerical solution, the toolbox uses the matrix $\exp\left(-\mathbf{A}\tau\right)$ as an initial estimation of $\mathbf{Q}_k$. In this case, that value turns out not to be in the region of attraction of a solution of \eqref{eq:findQ}.

\begin{figure}[!tb]
\centering
\includegraphics[scale=0.55]{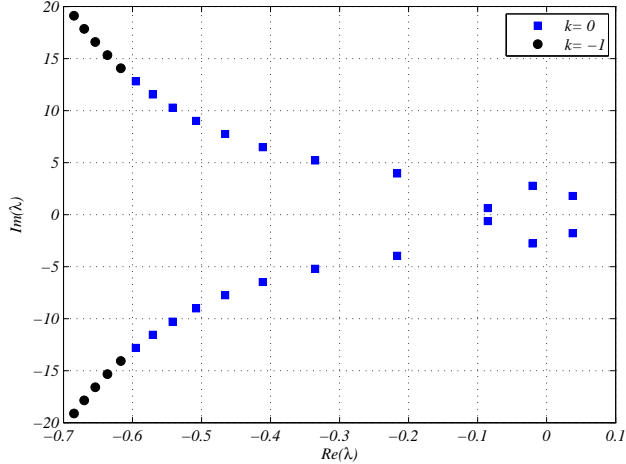}
\caption{Characteristic roots of the system under study. The roots represented as blue squares can be found using the principal branch of the matrix Lambert W function, whereas those in black circles are found using the branch corresponding to $k=-1$.}
\label{fig:roots}
\end{figure}

In order to obtain an \emph{a posteriori} guess for $\mathbf{Q}_k$, we reverse-engineered the solution, as in the proof of Proposition~\ref{prop1}. First, we use the QPmR algorithm \cite{Vhylidal2009} to find the characteristic roots of the system in a region close to the origin of the complex plane. The roots found are shown in Figure \ref{fig:roots}. 

The dominant roots of the system are $\lambda=0.0377\pm j1.7911$. Corresponding to these roots, we create the following $\mathbf{S}$ matrix:
\begin{equation}
\mathbf{S}=\left[\begin{array}{cc}0&1\\-3.2096&0.0753\end{array}\right]
\label{eq:s0}
\end{equation}
From \eqref{eq:s}, we have that
\begin{equation}
\mathbf{W}_k\left(\mathbf{M}\right)=\tau(\mathbf{S}-\mathbf{A})=\left[\begin{array}{cc}0&0\\8.9521&5.3766\end{array}\right].
\label{eq:w0}
\end{equation}
This shows that $W\left(m_{22}\right)\in\left[-1,\,\infty\right)$, which is the range of the principal branch of the Lambert W function. There is, therefore, a matrix $\mathbf{M}$ for which \eqref{eq:w0} is satisfied for $k=0$, and that matrix is
\begin{equation}
\mathbf{M}_0=\left[\begin{array}{cc}0&0\\1.9361&1.1628\end{array}\right]\times10^3.
\label{eq:M0}
\end{equation}
Since $\mathbf{B}$ and $\mathbf{M}$ are singular, there are an infinite number of $\mathbf{Q}_0$ matrices that satisfy $\mathbf{M}_0=\tau\mathbf{B}\mathbf{Q}_0$, for \eqref{eq:M0}. One of such matrices is
\begin{equation}
\mathbf{Q}_0=\left[\begin{array}{cc}1&1\\-650.3812&-392.6121\end{array}\right]
\label{eq:Q0}
\end{equation}
When this $\mathbf{Q}_0$ is used as starting value in the \emph{LambertDDE} toolbox for $k=0$, the numerical solution of \eqref{eq:findQ}, corresponding to the dominant roots, is found at the first iteration, as expected. Furthermore, if the matrix is slightly perturbed, the method still converges to the same solution after a few iterations.

Now, let us consider a non-dominant pair of roots: $\lambda=-0.4113\pm j6.4803$. Following a similar reasoning, we obtain the following matrices:
\begin{equation}
\begin{split}
\mathbf{S}&=\left[\begin{array}{cc}0&1\\-42.1633&-0.8226\end{array}\right],\\
\mathbf{W}_k\left(\mathbf{M}\right)&=\left[\begin{array}{cc}0&0\\-185.8166&0.8868\end{array}\right].
\end{split}
\label{eq:nondom}
\end{equation}
For this case, we again have $W\left(m_{22}\right)\in\left[-1,\,\infty\right)$. This implies that this pair of roots can also be found using the principal branch of the matrix Lambert W function. Indeed, if a matrix close to this one,
\begin{equation}
\mathbf{Q}_0=\left[\begin{array}{cc}1&1\\145.3412&-5.7175\end{array}\right],
\label{eq:Q01}
\end{equation}
which was created as in the previous case, is used as initial condition, the numerical routine within the \emph{LambertDDE} toolbox converges to the solution \emph{using the principal branch}.

In fact, we have observed that using suitable initial conditions, the 11 pairs of roots presented as blue squares in Figure \ref{fig:roots} can be found using the principal branch of the matrix Lambert W function. 

If we consider now the pair of roots $\lambda=-0.6169\pm j14.0734$, the corresponding $\mathbf{S}$ and $\mathbf{W}_k\left(\mathbf{M}_k\right)$ are
\begin{equation}
\begin{split}
\mathbf{S}&=\left[\begin{array}{cc}0&1\\-198.4405&-1.2338\end{array}\right],\\
\mathbf{W}_k\left(\mathbf{M}\right)&=\left[\begin{array}{cc}0&0\\-967.2027&-1.1692\end{array}\right].
\end{split}
\label{eq:km1}
\end{equation}
In this case, $W\left(m_{22}\right)\in\left(-\infty,-1\right]$, which is the range of the branch indexed by $k=-1$. Using this branch, we obtain a $\mathbf{Q}_{-1}$ matrix,
\begin{equation}
\mathbf{Q}_{-1}=\left[\begin{array}{cc}1&1\\95.1384&-4.8789\end{array}\right],
\label{eq:Qm1}
\end{equation}
which is a solution to \eqref{eq:findQ} for $k=-1$. Selecting an initial condition close to this matrix guarantees convergence to this solution. 

This procedure can be repeated for all the roots marked as black squares in Figure \ref{fig:roots}, as well as for roots further to the left of the complex plane, always using the branch corresponding to $k=-1$. This example shows how the whole spectrum of a system with the structure given in \eqref{eq:sys2nd} can be calculated using only two branches of the matrix Lambert W function and properly selected initial conditions for the solution of the nonlinear equation \eqref{eq:findQ}.

To show how higher branch numbers can also be used to find eigenvalues, let us consider a non conjugate eigenvalue pair, with $\lambda_1=-0.0204+j2.7705$ and $\lambda_2=-0.4658+j7.7500$. Notice how one of these eigenvalues was found using $k=0$, whereas the other was found using $k=-1$ in the previous exercise. With this eigenvalues the following matrix is created:
\begin{equation}
\mathbf{W}_k\left(\mathbf{M}\right)=\left[\begin{array}{cc}0&0\\132.3092+j07.2411&2.5693+j52.6026\end{array}\right]
\end{equation}
which has $W_k\left(m_{22}\right)$ in the range of the 9-th branch of the Lambert W function \cite{Corless1996}. Therefore this pair of eigenvalues can be found using $k=9$ and an appropriate initial condition in \eqref{eq:findQ}. Additionally, we have observed that keeping the same $\lambda_2$ and using the complex conjugate of $\lambda_1$, the matrix created is in the range of the branch indexed by $k=4$. This emphasizes the lack of an structured correspondence between the eigenvalues of \eqref{eq:sys2nd} and the branches of the Matrix Lambert W function.

\section{Odd Number of Real Characteristic Roots}
\label{sec:odd}
The previous discussions considered systems for which all the characteristic roots can be paired in such a way that a real $\mathbf{S}$ matrix is produced. This is the case when the complex roots come in conjugate pairs and the number of roots on the real axis is even. When there is an odd number on real roots, this is not possible. In this section, we undertake further discussions on this topic.
 
\begin{figure}[!tb]
\centering
\includegraphics[scale=0.55]{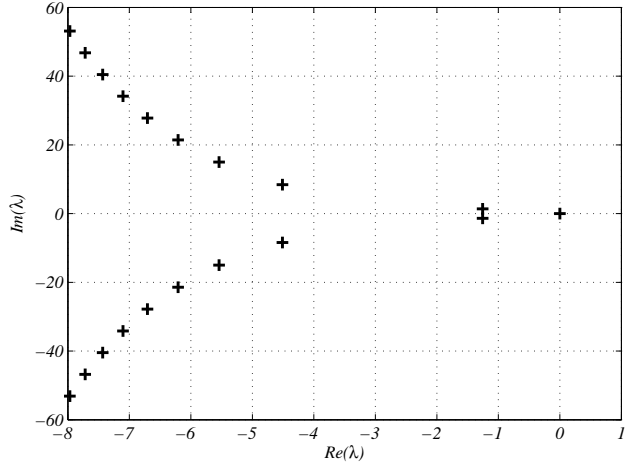}
\caption{Characteristic roots of the system in\eqref{eq:sys2}}
\label{fig:roots2}
\end{figure}

Consider the system with
\begin{equation}
\mathbf{A}=\left[\begin{array}{rc}0&1\\-1&0\end{array}\right],\ 
\mathbf{B}=\left[\begin{array}{cc}0&0\\1&0\end{array}\right],\ 
\tau=1,
\label{eq:sys2}
\end{equation}
whose characteristic equation is given by
\begin{equation}\label{eq:charvb}
\lambda^2+1-e^{-\lambda}=0.
\end{equation}
From \eqref{eq:charvb} it can be seen that there is only one real characteristic root, at the origin and with multiplicity one. Furthermore, this is the \emph{rightmost} root, as can be seen in Figure~\ref{fig:roots2}.

For this particular system $\mathbf{B}$ has a repeated zero eigenvalue. According to the observations in \cite{Yi2009,Yi2010d}, in this case the dominant root should be found using the principal branch or the branches $k=\pm1$. 

Using the ``reverse-engineering" approach for this example, we observe that all the complex conjugate eigenvalues can be obtained using $k=-1$ and a suitable initial condition. However, the dominant root at the origin cannot be found using the principal branch of the Lambert W function. The logic behind Propositions~\ref{prop1}-\ref{prop2} is based on constructing real $\mathbf{S}_k$, which cannot be done with this single real eigenvalue. A formal proof of this statement is given by contradiction. Assume that there is a real solution $\mathbf{M}_k$, and corresponding real $\mathbf{S}_k$, defined by (\ref{eq:5}), having eigenvalue $\lambda=0$. By (\ref{eq:5}) and (\ref{eq:findQ}) this implies that $\mathbf{S}_k$ satisfies (\ref{eq:char}). This implies on its turn that $(I,\mathbf{S})$ is an invariant pair of (\ref{eq:ltitds}), i.e.,~the eigenspace corresponding to the real characteristic roots is two-dimensional.

If we relax the conditions and allow $\mathbf{S}$ to be complex in the method presented earlier, we can match the characteristic root at the origin with any complex eigenvalue of the form $\lambda=a+jb$. This leads to matrices of the form:
\begin{align}
\mathbf{S}&=\left[\begin{array}{cc}0&1\\0&-a-jb\end{array}\right]\\
\mathbf{W}\left(\mathbf{M}\right)&=\left[\begin{array}{cc}0&0\\\tau&-\tau\left(a+jb\right)\end{array}\right]\label{eq:wsingle}.
\end{align}
The matrix given in \eqref{eq:wsingle} can be found using a value of $k$ such that $-j\tau b$ belongs to the range of the $k-$th branch of the Lambert W function. Therefore, for this particular system, the dominant root can be found using \emph{any} branch of the matrix Lambert W function, if a suitable initial condition is provided.

\section{Concluding Remarks} 
\label{sec:conclusions}
The methodology for the stability analysis and control synthesis of \textsc{lti-tds} presented in \cite{Yi2009,Yi2010d}, based on the matrix Lambert W function, assumes that there is a one to one correspondence between the characteristic roots of the system and the branches of the matrix Lambert W function. It also assumes that this correspondence is such that using the principal branch always leads to the dominant roots of the system. This paper illustrates that such a correspondence does not exist in all the cases, and that for a particular, albeit very common structure for the system, the branches corresponding to $k=0$ and $k=-1$ can be used to find all the characteristic roots of the system. It is also illustrated that the correspondence cannot always be restored when selecting the particular initial condition $\exp(-\mathbf{A}\tau)$ in solving \eqref{eq:findQ}.

In our examples matrix $\mathbf{B}$ had reduced rank. An interesting path of research is to find out whether the correspondence between characteristic roots and branches of the Lambert W function, holding for scalar systems, can be extended to higher-order systems, provided \emph{additional} structural conditions on the system are assumed.

The MATLAB code used to create the examples can be requested via email to the authors or can be downloaded from \texttt{https://db.tt/mSI3VwbO}.

\section*{Acknowledgments}
This paper is the result of a research visit that R. Cepeda-Gomez made to the Computer Science department at KU Leuven, funded by a Coimbra Group Scholarship for Young Professors and Researchers of Latin America. This scholarship was awarded while he was a lecturer of the Mechatronic Engineering Department at Universidad Santo Tom\'as in Bucaramanga, Colombia.

The work of W. Michiels work is supported by the Programme of Interuniversity Attraction Poles of the Belgian
Federal Science Policy Office (IAP P6-DYSCO), by OPTEC, the Optimization in Engineering Center
of KU Leuven, and by the project G.0712.11N of the Research Foundation-Flanders (FWO).

The authors thank the comments received from the anonymous reviewers, which helped to improve the presentation of the paper, specially in the second example.

\bibliographystyle{plain}
\bibliography{lambert}
\end{document}